\documentclass{article}
\usepackage{graphicx}
\usepackage{amsfonts}
\usepackage{amsmath}
\usepackage{amsmath, amscd, amsfonts, amsthm, tikz, times}

\newcounter{conjecture}
\newcounter{remark}
\newcounter{corollary}

\newenvironment{corollary}{\medskip{\bf Corollary \thecorollary.}
\addtocounter{corollary}{1}\em}{\rm}
\newtheorem{theorem}{Theorem}
\newtheorem{lemma}{Lemma}

\newcommand{\dd}{\delta}

\newcommand{\lar}{\longrightarrow}
\newcommand{\eps}{\varepsilon}

\newcommand{\reals}{\mathbb{R}}
\newcommand {\rrr}[1]{(\ref{#1})}

\def \be{\begin{equation}}
\def \ee{\end{equation}}
\def \bt{\begin{theorem}}
\def \et{\end{theorem}}
\def \bc{\begin{corollary}}
\def \ec{\end{corollary}}
\def \bea{\begin{eqnarray}}
\def \eea{\end{eqnarray}}
\def \bas{\begin{eqnarray*}}
\def \eas{\end{eqnarray*}}
\def \bl{\begin{lemma}}
\def \el{\end{lemma}}

\def \ga{\gamma}

\def \Om{\Omega}

\def \vski{\vspace{12pt}}
\def \ff{\infty}

\def \AA{{\cal A}}

\def \DD{\Delta}

\def \RR{{\mathbb R}}

\def \reals{{\mathbb R}}

\def \sgn{\mathop{\mathrm{sgn}}}

\def \({\left(}
\def \){\right)}

\def \bc{\begin{center} }
\def \ec{\end{center} }
\def \bs{\begin{slide} }
\def \es{\end{slide} }

\def\square{{\vcenter{\vbox{\hrule height.3pt
         \hbox{\vrule width.3pt height5pt \kern5pt
            \vrule width.3pt}
         \hrule height.3pt}}}}
\def\qed{{\hfill $\square$ \bigskip}}

\newcounter{cccases}
\begin{document}

\title{Symmetry in the Green's function for birth-death chains}
\author{Greg Markowsky\\School of Mathematical Sciences\\Monash University, Melbourne, Australia\\and\\Jos\'e Luis Palacios\\ Department of Electrical and Computer Engineering\\  The University of New Mexico, Albuquerque, USA.}



\maketitle

\begin{abstract}
A symmetric relation in the probabilistic Green's function for birth-death chains is explored. Two proofs are given, each of which makes use of the known symmetry of the Green's functions in other contexts. The first uses as primary tool the local time of Brownian motion, while the second uses the reciprocity principle from electric network theory. We also show that the the second proof extends easily to cover birth-death chains (a.k.a. state-dependent random walks) on trees.
\end{abstract}


\section{Introduction}
One of the most intriguing properties of the classical Green's function $G_\Om$ for a domain $\Om$ in $\RR^n$ is symmetry: $G_\Om(x,y) =G_\Om(y,x)$. This was originally noted in the physical setting, but in relatively more recent years the moniker "Green's function" has been adopted by probabilists in reference to the expected amount of time a random process spends at a point. This is because, under appropriate conditions, the Green's function in $\RR^n$ can be realized as

\begin{equation} \label{yi}
G_\Om(x,y) = \int_{0}^{\ff} \rho_t(x,y)dt,
\end{equation}

where $\rho_t(x,y)$ is the density at point $y$ at time $t$ of a Brownian motion starting at $x$ and killed upon exiting $\Om$. More generally, in probabilistic settings the Green's function is generally taken to be density of the occupation measure

\begin{equation} \label{yi2}
\mu_x(A) = E_x \int_{0}^{\ff} 1_A (X_t)dt
\end{equation}

when $X_t$ is a continuous-time stochastic process (which is equivalent to \rrr{yi} when $X_t$ is Brownian motion), or

\begin{equation} \label{yi3}
G(x,y) = E_x\sum_{m=0}^{\ff} 1_{\{X_m=y\}}
\end{equation}

when $X_m$ is a discrete-time random process. It is difficult to find an intutive probabilistic reason for the symmetry property exhibited by $G_\Om$ defined by \rrr{yi} for general domains in $\RR^n$, and such a symmetry does not hold for the Green's functions \rrr{yi2} and \rrr{yi3} defined by more general processes. In this paper, we are interested in investigating the symmetry property for a general class of discrete-time processes known as {\it birth-death chains}. \\

A birth-death chain $X_m$ is a Markov chain taking values on the integers with the following transition probabilities:

\be p_{nj} = \left \{ \begin{array}{ll}
r_{n} & \qquad  \mbox{if } j=n+1  \\
l_n & \qquad \mbox{if } j=n-1 \\
a_{n} & \qquad  \mbox{if } j=n  \\
0 & \qquad \mbox{if } |n-j| > 1\;,
\end{array} \right. \ee

with $l_n+r_n+a_n=1$. In order to avoid a number of qualifying statements attached to our results, we assume that $r_n, l_n >0$ for all $n$; the methods applied below can easily be adapted to any special case under investigation in which this does not hold. We allow the chain $X_m$ to have absorbing states if desired, so that there may be values $A$ and/or $B$ such that $X_m$ is killed upon reaching $A$ or $B$; to be precise, we set $X_m = \DD$, where $\DD$ is a cemetery point, for all $m \geq \tau = \inf \{m':X_{m'} = A \mbox{ or } B\}$, with $\tau$ taken to be $\ff$ on the set $\{X_{m'} \neq a \mbox{ or } b \mbox{ for all } m' \geq 0\}$ or in the case in which absorbing states are not present. Again to reduce the number of qualifying statements, we assume that if there are two such absorbing states then for any initial points $j,k$ of the walk we have $A < j,k < B$, or if there is only one absorbing state $B$ then either $j,k <B$ or $j,k>B$, etc. We define the Green's function $G(x,y)$ as in \rrr{yi3} by

\begin{equation} \label{yi4}
G(x,y) = E_x\sum_{m=0}^{\tau} 1_{\{X_m=y\}},
\end{equation}

the expected number of visits to $y$ of the chain which starts at $x$ before time $\tau$. We will prove the following theorem.

\begin{theorem} \label{ima} Assume $j<k$. If the chain $X_m$ is recurrent, then $G(j,k) = G(k,j) = \ff$, but if $X_m$ is transient then

\begin{equation} \label{kyoto}
\begin{split}
\frac{G(j,k)}{G(k,j)} = \frac{r_{j} \ldots r_{k-1}}{l_{j+1} \ldots l_{k}}.
\end{split}
\end{equation}

\end{theorem}

{\bf Remark:} It is interesting to note that the relation \rrr{kyoto} is independent of the behavior of the chain at states $n$ with $n < j$ or $n > k$, and holds whether or not absorbing states are present. It may also seem to be the case (and be hard to believe) that the quantities $a_n$ should play no role here; however this is a bit illusory since the relation between $a_j$ and $a_k$ is in essence carried by the ratio $\frac{r_j}{l_k}$ in \rrr{kyoto}, as is shown at the beginning of the next section.

\vski

We will give two proofs of Theorem \ref{ima}. The first makes use of Brownian motion and its corresponding theory of local time, and in the end depends upon the symmetry of the classical Green's function for an interval in $\RR$. The second makes use of electric network theory, and makes use of the symmetry property of voltages known as the reciprocity principle (which can be translated to the symmetry of the classical Green's function when voltages are taken in domains in $\RR^n$). The next section contains these two proofs. The electric-networks proof extends easily to the case of birth-death chains on trees, and the final section contains the necessary details on this, as well as a few examples.

\section{Proofs of Theorem \ref{ima}.} \label{cha}

Both proofs which are to be presented deal most naturally with the situation $a_n = 0$ for all $n$, so we begin by showing why we may assume this. Let us suppose that the result holds whenever $a_n = 0$ for all $n$, and let $X_m$ be an arbitrary chain (where the $a_n$'s are allowed to be nonzero). Define a new birth-death chain $\tilde X_m$ with transition probabilities $\tilde a_n = 0, \tilde l_n = \frac{l_n}{l_n+r_n}, \tilde r_n = \frac{r_n}{l_n+r_n}$. $\tilde X_m$ is simply the chain $X_m$ but with the waiting times between moves removed, and it follows that if $\tilde G$ is the Green's function of $\tilde X_m$ then $G(j,k) = \frac{1}{1-a_k}\tilde G(j,k)$; this is a consequence of the fact that the expected time until first success of a Bernoulli trial with probability of success $p$ is $\frac{1}{p}$, so that the expected time for $X_m$ to make a nontrivial move upon reaching state $k$ is $\frac{1}{1-a_k}$. We then have

\begin{equation} \label{}
\begin{split}
\frac{G(j,k)}{G(k,j)} &= \frac{(1-a_j)\tilde G(j,k)}{(1-a_k) \tilde G(k,j)} = \Big(\frac{1-a_j}{1-a_k} \Big) \Big( \frac{\tilde r_j}{\tilde l_k}\Big) \Big( \frac{\tilde r_{j+1} \ldots \tilde r_{k-1}}{\tilde l_{j+1} \ldots \tilde l_{k-1}} \Big) \\
& = \Big(\frac{1-a_j}{1-a_k} \Big) \Big( \frac{r_j/(r_j+l_j)}{l_k/(r_k+l_k)}\Big) \Big( \frac{r_{j+1} \ldots r_{k-1}}{l_{j+1} \ldots l_{k-1}}\Big) = \frac{r_{j} \ldots r_{k-1}}{l_{j+1} \ldots l_{k}},
\end{split}
\end{equation}

where we have used $1-a_n = r_n+l_n$ and $\frac{\tilde r_n}{\tilde l_n} = \frac{r_n}{l_n}$. We may therefore assume $a_n = 0$ for all $n$ in what follows, and begin the proofs proper.

\vski

{\it Brownian motion proof:} In \cite{me} and \cite{me2}, a technique was given for realizing $X_m$ as a Brownian motion stopped at a certain sequence $\tau(m)$ of stopping times, and we now briefly discuss (a slightly modified version of) the technique. Set $x_0 = 0$, $t_0 = l_0$, and for $n \geq 1$ set

\be \label{pred}
t_n := \frac{l_0 l_1 l_2 \ldots l_n}{r_1 r_2 \ldots r_n}
\ee

and $x_n = \sum_{j=0}^{n-1} t_j$. Set $t_{-1} = -r_0$, for $n \leq -2$ set

\be \label{pred}
t_n := -\frac{r_0 r_{-1} r_{-2} \ldots r_{n+1}}{l_{-1} l_{-2} \ldots l_{n+1}},
\ee

and for $n \leq -1$ set $x_n = \sum_{j=n}^{-1} t_j$ (note that $x_n < 0$ when $n<0$). Since the sequence $\{x_n\}_{n=-\ff}^\ff$ is increasing in $n$ it converges on either side to upper and lower limits $x^-_{\ff}$ and $x^+_\ff$, one or both of which may be infinite. Let $B_t$ be a Brownian motion stopped at the first time $\tau(\Delta)$ it hits $x^-_\ff$ or $x^+_\ff$. Set ${\cal A} = \cup_{n=-\ff}^\ff \{x_n\}$. We will be starting the Brownian motion at a point in $\AA$, and define the stopping times $\tau(m)$ recursively by setting $\tau(0)=0$, and having defined $\tau(m)$ we let $\tau(m+1)=\inf_{t>\tau(m)}\{B_t \in \cal{A},$ $B_t \neq B_{\tau(m)}\}$. That is, the $\tau(m)$'s are the the successive hitting times of points in $\cal{A}$. We see that the variables $B_{\tau(0)},B_{\tau(1)},B_{\tau(2)}, \ldots $ form a random process taking values in $\cal{A}$. Let $\phi:\cal{A} \lar \reals$ be defined by $\phi(x_n)=n$. The strong Markov property of Brownian motion and the formula for the exit distribution of Brownian motion from an interval imply that $\phi(B_{\tau(0)}), \phi(B_{\tau(1)}), \phi(B_{\tau(2)}), \ldots$ is a realization of our birth-death chain (\cite{me}). We may therefore take $X_m = \phi(B_{\tau(m)})$ in what follows. An important quantity for us will be the {\it local time} of Brownian motion, which is the density of the occupation measure of Brownian motion with respect to Lebesgue measure. That is, the local time $L_t^x$ satisfies

\be \label{}
L_t^x dx = \int_{0}^{t} 1_{B_s \in dx} ds
\ee

It is well known that $L_t^x$ exists and that

\be \label{bond}
L_t^x = \lim_{\eps \lar 0} \frac{1}{2\eps} \int_{0}^{t} 1_{|B_s-x|<\eps}ds
\ee

almost surely. Formally, we write

\be \label{}
L_t^x = \int_{0}^{t} \dd_x(B_s)ds
\ee

in place of \rrr{bond}, where $\dd_x$ is the Dirac delta function. The local time provides a measure of the amount of time that Brownian motion spends at a point, and perhaps not surprisingly we can make use of it in our study of the number of visits that a birth-death chain makes to a given state. The following lemma appears in \cite{me2}, but we include the short proof for the benefit of the reader.

\begin{lemma} \label{rer}

Suppose $a< y,z < b$. Let $\ga(a,b) = \inf_{t>0} \{B_t =a \mbox{ or }b\}$. Then

\begin{equation} \label{wwd}
E_z[L^y_{\ga(a,b)}] = \frac{(b-z)(y-a)+(z-a)(b-y)}{b-a} - |z-y|.
\end{equation}
\end{lemma}

{\bf Proof:} Tanaka's formula (see \cite{fima} or \cite{rosmarc}) states that

\begin{equation} \label{}
|B_t - y| = |B_0 - y| + \int_{0}^{t}\sgn(B_s-y) dB_s + L_t^y.
\end{equation}

The stochastic integral here is a martingale, and the optional stopping theorem can be applied. We find

\begin{equation} \label{}
E_z[L^y_{\ga(a,b)}] = E_z|B_{\ga(a,b)} - y| - |z-y|
\end{equation}

Since $P_z(B_{\ga(a,b)} = a) = \frac{b-z}{b-a}$ and $P_z(B_{\ga(a,b)} = b) = \frac{z-a}{b-a}$, we see that $E_z|B_{\ga(a,b)} - y| = \frac{b-z}{b-a}(y-a) + \frac{z-a}{b-a}(b-y)$, and the result follows. \qed

The formula \rrr{wwd} presents several features of interest. First, we note that $E_z[L^y_{\ga(a,b)}] = E_y[L^z_{\ga(a,b)}]$, and this is in fact a special case of the symmetry of the classical Green's function, as the quantity $E_y[L^z_{\ga(a,b)}]$ is equal to $G_{(a,b)}(y,z)$ defined in \rrr{yi}. To see this formally, use the identity $E_y[\phi(B_t)] = \int_{-\ff}^{\ff}\phi(x)\rho_t(y,x)dx$, and write

\begin{equation} \label{}
E_y[L^z_{\ga(a,b)}] = E_y \int_{0}^{\ff} \dd_z(B_t)dt = \int_{-\ff}^{\ff}\Big(\int_{0}^{\ff}\dd_z(x) dt \Big)\rho_t(y,x)  dx = \int_{0}^{\ff} \rho_t(y,z) dt,
\end{equation}

where $B_t$ is Brownian motion killed upon leaving $(a,b)$, and $\rho_t$ is its corresponding density (these manipulations can be made rigorous with little difficulty). The second noteworthy feature of \rrr{wwd} is that, when $z$ and $y$ coincide, the formula \rrr{wwd} simplifies substantially to $E_z[L^z_{\ga(a,b)}] = \frac{2(b-z)(z-a)}{b-a}$. This expression will be important in what follows.

\vski

Now, if we start our Brownian motion at $x_j$ then the total expected local time accumulated at $x_k$ before hitting $x^-_\ff$ or $x^+_\ff$ is given by $E_{x_j}[L_{\ga(x^-_\ff,x^+_\ff)}^{x_k}]$. However, if we condition on $B_{\tau(m)} = x_k$ we see that the expected local time accumulated at $x_k$ between times $\tau(m)$ and $\tau(m+1)$ is $E_{x_k}[L_{\ga(x_{k-1},x_{k+1})}^{x_k}]$. This is the amount of local time accumulated by $B_t$ at $x_k$ each time the birth-death chain $\phi(B_t)$ visits the point $k$. We therefore obtain

\begin{equation} \label{edge}
G(j,k) = \frac{E_{x_j}[L_{\ga(x^-_\ff,x^+_\ff)}^{x_k}]}{E_{x_k}[L_{\ga(x_{k-1},x_{k+1})}^{x_k}]}.
\end{equation}

Using $x_n - x_{n-1} = t_{n-1}$ and Lemma \ref{rer} we obtain

\begin{equation} \label{}
\begin{split}
G(j,k) \frac{2t_{k-1}t_k}{t_{k-1} + t_k} & = G(j,k)E_{x_k}[L_{\ga(x_{k-1},x_{k+1})}^{x_k}] = E_{x_j}[L_{\ga(x^-_\ff,x^+_\ff)}^{x_k}] = E_{x_k}[L_{\ga(x^-_\ff,x^+_\ff)}^{x_j}] \\
& = G(k,j)E_{x_j}[L_{\ga(x_{k-1},x_{k+1})}^{x_j}] = G(k,j) \frac{2t_{j-1}t_j}{t_{j-1} + t_j}.
\end{split}
\end{equation}

In order to simplify we will use the identity $\frac{t_n}{t_{n'}} = \frac{l_{n'+1}l_{n'+2}\ldots l_{n-1}l_n}{r_{n'+1}r_{n'+2}\ldots r_{n-1}r_n}$, valid whenever $n'<n$. We obtain

\begin{equation} \label{}
\begin{split}
\frac{G(j,k)}{G(k,j)} &= \Big(\frac{t_{k-1} + t_k}{t_{j-1} + t_j} \Big)\Big(\frac{t_{j-1}t_j}{t_{k-1}t_k}\Big) \\
& = \Big(\frac{\frac{l_j \ldots l_{k-1}}{r_j \ldots r_{k-1}} + \frac{l_j \ldots l_{k}}{r_j \ldots r_{k}}}{1 + \frac{l_j}{r_j}}\Big)\Big( \frac{r_j \ldots r_{k-1}}{l_j \ldots l_{k-1}} \Big)\Big(\frac{r_{j+1} \ldots r_{k}}{l_{j+1} \ldots l_{k}}\Big) \\
& = \Big(\frac{r_k+l_k}{r_j + l_j}\Big) \Big(\frac{r_{j} \ldots r_{k-1}}{l_{j+1} \ldots l_{k}}\Big) = \frac{r_{j} \ldots r_{k-1}}{l_{j+1} \ldots l_{k}},
\end{split}
\end{equation}

where in the second equality we have divided top and bottom by $t_{j-1}$ and in the final one we have used $r_k+l_k = r_j + l_j =1$. This completes the proof. $\bullet$
\vskip .2 in
{\it Electric preliminaries.}
\vskip .2 in
  On a finite connected undirected graph $G=(V,E)$ with $|G|=n$ such that the edge between
vertices $i$ and $j$ is given a
resistance $r_{ij}$ (or equivalently, a conductance $C_{ij}=1/r_{ij}$),
we can define
 the random walk on $G$ as the Markov chain  $X_n, n \ge 0$, that from
its current vertex $v$ jumps
to the neighboring vertex $w$ with probability  $p_{vw}=C_{vw}/C(v)$,
where $C(v)=\sum_{w: w\sim v} C_{vw}$, and $w\sim v$ means that $w$ is a
neighbor of $v$.  There may be a
conductance $C_{zz}$ from a vertex $z$
to itself, giving rise to a transition probability form z to itself (such as $a_z$ from the previous section), though the most studied case of these random walks, the simple random walk, excludes the loops and considers all $r_{ij}$'s to be equal to 1.

\vski

The beginner's handbook when studying random walks on graphs from the
viewpoint of electric networks is \cite{doysne},  which is both a mandatory reference in this area of research and a textbook suitable for undergraduate students.  We begin with a basic fact from that text which we state as a lemma.  Consider a general random walk on a finite graph,  and let $E_aN_x^b$ be the expected number of times the vertex $x$ is visited by the walk started at $a$ before it reaches $b$.  Then we have

\begin{lemma}
\begin{equation}
\label{esc2}
E_aN_x^b=C(x)V_x,
\end{equation}
where $V_x$ is the voltage at $x$ when a battery is placed between $a$ and $b$ such that the current entering at $a$ is 1 and the voltage at $b$ is 0.
\end{lemma}

A crucial electrical result we need below is the reciprocity principle, a simple proof of which  (based on the superposition principle for electric networks) can be found in \cite{palhit}.
\begin{lemma}
For any $a, b, c \in G$, if we set a battery between $a$ and $b$ so that a unit current enters $a$ and exits at $b$ and the voltage at $b$ is $0$, then the value of the voltage at $c$ is the same as the value of the voltage at $a$ when we disconnect the battery cable at $a$ and reconnect it at $c$.
\end{lemma}

{\it Electric Proof of the theorem}. In \cite{paltet} it was shown that every birth-and-death process on the integers can be expressed as a random walk on the linear graph endowed with the conductances
$C_{0,0}=a_0, C_{0,1}=r_0, C_{-1,0}=l_0$,


$$ C_{zz}= \left \{ \begin{array}{ll}
a_z\frac{r_0r_1\cdots r_{z-1}}{l_1l_2\cdots l_z}{\rm ~if~} z \geq 1 \\
a_z\frac{l_0l_{-1}\cdots l_{z+1}} {r_{-1}r_{-2}\cdots r_z} & {\rm ~if~} z \leq -1\;, \end{array} \right. $$
and
$$ C_{z,z+1}= \left \{ \begin{array}{ll}
\frac{r_0r_1\cdots r_z}{l_1l_2\cdots l_z} & {\rm ~if~} z\ge 1 \\
\frac{l_0l_{-1}\cdots l_{z+1}} {r_{-1}r_{-2}\cdots r_{z+1}} & {\rm ~if~} z\le -2.\end{array} \right.$$

Some algebra yields
$$C(z)=\left \{ \begin{array}{ll}
\frac{r_0\cdots r_{z-1}}{l_1\cdots l_z} & {\rm ~if~} z\ge 1 \\
1 & {\rm ~if~} z=0 \\
\frac{l_0\cdots l_{z+1}}{r_{-1}\cdots r_z} & {\rm ~if~} z \le -1.\end{array} \right.$$

Now if for the birth-and-death process we look at the problem of finding   $\displaystyle G(x,y)=E_x\sum_{m=0}^\tau 1_{\{X_m=y\}}$, where $\tau=\inf \{n: X_m \in \{a,b\}\}$, and $a<x<y<b$, by lemma 2 this is equivalent to finding the product $C(y)V(y)$, where $V(y)$ is the voltage at $y$ when we inject a unit current at $x$ and establish a zero voltage at both vertices $a$ and $b$, and where $C(y)$ is the total conductance emanating from $y$.   By the reciprocity principle we have
$$V(y)=V^{'}(x),$$
where $V^{'}(x)$ is the voltage at $x$ when we inject a unit current at $y$ and keep a zero voltage al $a$ and $b$.  Or equivalently,
\begin{equation}
\label{previous}
\frac{G(x,y)}{G(y,x)}=\frac{C(y)}{C(x)}.
\end{equation}
Now, the closed form expressions for $C(z)$ given above allow us to express (\ref{previous}) as
$$\frac{C(y)}{C(x)}=\frac{r_x\cdots r_{y-1}}{l_{x+1}\cdots l_y},$$
thus proving (\ref{kyoto}). $\bullet$
\vskip .2 in
It should be noted that (\ref{previous})  holds for a general graph endowed with arbitrary conductances, and this will allow us to extend our results to birth-death chains on arbitrary trees in the next section.

\section{Green's functions for birth-death chains on trees}

Consider a birth-death chain on a tree, that is, a Markov chain where the transitions occur from a given vertex of a (possibly infinite) tree either to itself or to any neighboring vertex in the tree.  Consider a {\it finite} subtree $T$ of the given tree and let $I$ and $L$ be the sets of interior points and leaves, respectively, of $T$.   Take $x, y\in I$ and define $G(x,y)$ as in (\ref{yi4}) except that now we take
$$\tau=\inf\{m: X_m\in L\}.$$

We look again at the problem of expressing the quotient
$\displaystyle \frac{G(j,k)}{G(k,j)}$ in terms of the transition probabilities, as in (\ref{kyoto}), and for that purpose we will use (\ref{previous}) and the fact that any birth-death chain on a finite tree is expressible as a random walk on the underlying tree endowed with appropriate conductances found with a simple algorithm essentially equal to the one found in \cite{palqir}, that we present here for completeness:
\vskip .2 in
1. Look at the subtree $T$ of the original tree $\hat{T}$, and delete all transition probabilities $p(i,j)$, where $i\in \hat{T}-T$, $j\in L$, as well as those where $i\in L, j\in T$.

2. Take any vertex $v \in V$ of the tree as the root, and consider any $u\sim v$. Assign an arbitrary (positive) value to $C_{vu}$.

3. Obtain the conductance $C_{vv}$ and all conductances $C_{vw}$ where $w$ is a neighbor of $v$ as

\begin{equation}
\label{Csub}
 C_{vw}={{C_{vu}p(v,w)}\over p(v,u)}, ~w\sim v;~~~~ C_{vv}={{C_{vu}p(v,v)}\over {p(v,u)}}.
 \end{equation}

4. Taking $v$ as the root, traverse the vertices of the tree using Breadth First Search (BFS). Every time a vertex not previously visited is reached, only one of its adjacent conductances has being assigned. Take this conductance as the $C_{vu}$ used in step 2 to obtain all other adjacent ones.

\vskip .2 in
The fact that the procedure works, that is, the fact that we can recover the transition probabilities from the conductances can be checked easily since
\begin{equation}
\label{C}
C(v)=\sum_{w} C_{vw}={{C_{vu}}\over {p(v,u)}},
\end{equation}
and then the motion of the random walk from $v$ to $w$ is dictated by
$${{C_{vw}}\over {C(v)}}={{C_{vu}p(v,w)}\over {p(v,u)}}{{p(v,u)}\over {C_{vu}}}=p(v,w),$$
when $v\sim w$ and
$${{C_{vv}}\over {C(v)}}={{C_{uv}p(v,v)}\over {p(v,u)}}{{p(v,u)}\over {C_{vu}}}=p(v,v),$$
as desired. This procedure stops when all leaves, and thus all vertices,  have been visited. Now we can prove the following generalization of formula (\ref{kyoto})
\begin{theorem} \label{fuchu}
Let the unique path between $j$ and $k$ be determined by the vertices $v_1, \ldots, v_k$, so that $j\sim v_1$, $v_i\sim v_{i+1}$,  for $1\le i \le k-1$,  and $v_k\sim k$.
Then
\begin{equation}
\label{final}
\frac{G(j,k)}{G(k,j)}=\frac{p(j,v_1)p(v_1,v_2)\cdots p(v_{k-1},v_k)p(v_k,k)}{p(k,v_k)p(v_k,v_{k-1})\cdots p(v_2,v_1)p(v_1,j)}
\end{equation}
\end{theorem}

{\bf Proof.} Assume first that $j\sim k$. We start the algorithm assigning $C_{jk}=1$.  Then, by  (\ref{C}) we have $C(j)=\frac{C_{jk}}{p(j,k)}=\frac{1}{p(j,k)}$ and $C(k)=\frac{C_{kj}}{p(k,j)}=\frac{1}{p(k,j)}$.  Then (\ref{previous}) implies
$$\frac{G(j,k)}{G(k,j)}=\frac{C(k)}{C(j)}=\frac{p(j,k)}{p(k,j)}.$$
Assume now that $j, v_1, k$ is the unique path between $j$ and $k$.  Again, we start the algorithm assigning $C_{jk}=1$.  Then, by  (\ref{C}) we have $C(j)=\frac{C_{jv_1}}{p(j,v_1)}=\frac{1}{p(j,v_1)}$. Now $C(k)=\frac{C_{kv_1}}{p(k,v_1)}=\frac{C_{v_1k}}{p(k,v_1)}$, and using (\ref{Csub}) in order to write $C_{v_1k}$ in terms of $C_{v_1j}$ we obtain
$$C(k)=\frac{C_{v_1j}p(v_1k)}{p(k,v_1)p(v_1,j)}=\frac{p(v_1,k)}{p(k,v_1)p(v_1,j)}.$$
Again by (\ref{previous}) we get
$$\frac{G(j,k)}{G(k,j)}=\frac{C(k)}{C(j)}=\frac{p(j,v_1)p(v_1,k)}{p(k,v_1)p(v_1,j)}.$$

It should be clear how to proceed by induction using (\ref{Csub}) and (\ref{C}). $\bullet$
\vskip .2 in
Notice that the value of $\frac{G(j,k)}{G(k,j)}$ does not depend explicitly on the number of branching points between $j$ and $k$, or the number of leaves of the tree, and its computation does not require to know the values of the conductances, only the transitions of going from $j$ to $k$ through the unique path and vice versa.
\vskip .2 in
{\bf Example:}  Assume the subtree $T$ in question is given in Figure 1. Then, for instance $\displaystyle \frac{G(2,4)}{G(4,2)}=\frac{p(2,3)p(3,4)}{p(4,3)p(3,2)}=\frac{36}{5}.$
\begin{figure}[h!] 
  \centering
  \includegraphics[width=5 in,height=2.5 in]{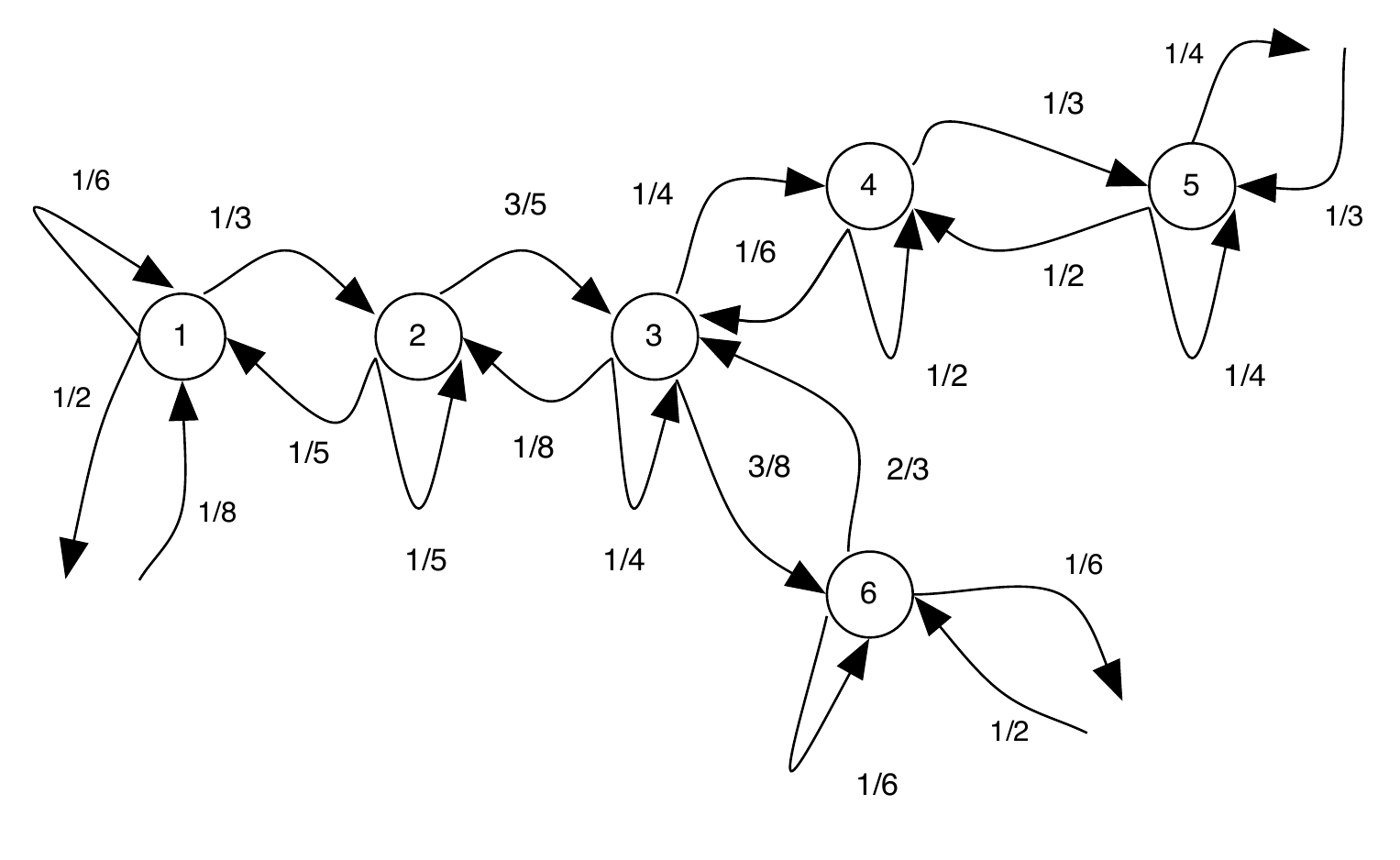}
\caption {The subtree}
\end{figure}

\vski

{\bf Remark:} The Brownian motion proof given in Section \ref{cha} does not seem to apply to Theorem \ref{fuchu}.

\vskip .2 in
\section{Acknowledgements}

The first author would like to express his gratitude for support from Australian Research Council Grants DP0988483 and DE140101201.

\bibliographystyle{alpha}
\bibliography{prob}

\end{document}